# Positive systems analysis via integral linear constraints

Sei Zhen Khong, Corentin Briat and Anders Rantzer

*Abstract*— Closed-loop positivity of feedback interconnections of positive monotone nonlinear systems is investigated. It is shown that an instantaneous gain condition on the open-loop systems which implies feedback well-posedness also guarantees feedback positivity. Furthermore, the notion of integral linear constraints (ILC) is utilised as a tool to characterise uncertainty in positive feedback systems. Robustness analysis of positive linear time-varying and nonlinear feedback systems is studied using ILC, paralleling the well-known results based on integral quadratic constraints.

## I. INTRODUCTION

Monotone systems are a class of dynamical systems that preserve a closed order relation from the input space and initial state to their state space and output space [1]. They constitute one of the most important classes used in mathematical biology and chemical modelling. Since biological and chemical models often treat variables such as population densities or concentrations of chemical mixtures that are intrinsically positive, they typically preserve positivity of solutions to differential equations and enjoy additional monotonicity or order-preserving properties.

Positive monotone systems have gained an increasing attention over the last decade due to the fact that their order-preserving properties can often be exploited to simplify computations and control synthesis. For instance, stabilising output feedback controllers for positive systems are characterised using linear programming in [2] and extensions to input-output gain optimisation are considered in [3]. In [4], it is shown that the input-output gain of positive systems can be evaluated using a diagonal quadratic storage function and this is utilised for $H_\infty$ optimisation of decentralised controllers in terms of semidefinite programming. Rantzer [5] establishes that positive systems are amenable to distributed control whose complexity scales linearly with the number of interconnections. Robust stability of positive monotone systems subject to time-varying delays is studied in [6]. The class of differentially positive systems is introduced in [7] as an extension of the definition of monotone systems.

Of the aforementioned works, [3] introduces the use of integral linear constraints (ILCs) to characterise uncertainty and derives linear-programming conditions to verify robust stability of LTI systems based on a dissipativity approach. This paper further generalises this idea via an input-output approach to analyse robust stability of positive feedback systems, in a similar spirit to integral quadratic constraints (IQCs) based analysis [8]. IQCs are a well-known and useful tool for describing uncertain input-output dynamics. It unifies a range of methods in the literature, such as small-gain and passivity. This paper presents a first attempt to inherit these features via the more scalable form of ILCs for positive feedback systems. The main result generalises the static gain conditions on robust stability of positive feedback interconnections of linear time-invariant (LTI) systems to linear time-varying (LTV) and nonlinear systems via integral linear constraints. Importantly, a sufficient condition guaranteeing positivity of the feedback interconnection of two positive nonlinear systems is also provided.

The paper evolves according to the following structure. Notation and preliminary material are introduced in the next section. Section III establishes a sufficient condition under which positive feedback interconnections preserve positivity of open-loop nonlinear systems. ILCs are reviewed in Section IV. Robust stability analysis via ILCs is performed in Section V. Section VI contains a couple of illustrative examples. Some conclusions are provided at the end.

## II. NOTATION AND PRELIMINARIES

### A. Signals and systems

Denote by $\mathbf{L}_1^n$ the set of $\mathbb{R}^n$-valued Lebesgue integral functions:
$$\mathbf{L}_1^n := \left\{ v : [0,\infty) \to \mathbb{R}^n : \|v\|_1 := \int_0^\infty |v(t)|\, dt < \infty \right\},$$
where $|\cdot|$ denotes the 1-norm. Let $\mathbb{R}_+^n$ denote the nonnegative orthant of $\mathbb{R}^n$. Define
$$\mathbf{L}_{1+}^n := \{v \in \mathbf{L}_1 : v(t) \in \mathbb{R}_+^n \text{ a.e.}\},$$
where "a.e." is with respect to the Lebesgue measure on $[0,\infty)$. Define the truncation operator
$$(P_T v)(t) := \begin{cases} v(t) & t \in [0,T) \\ 0 & \text{otherwise}, \end{cases}$$
and the extended spaces
$$\mathbf{L}_{1e}^n := \{v : [0,\infty) \to \mathbb{R}^n : P_T v \in \mathbf{L}_1 \ \forall T \in [0,\infty)\};$$
$$\mathbf{L}_{1e+}^n := \{v : [0,\infty) \to \mathbb{R}^n : P_T v \in \mathbf{L}_{1+} \ \forall T \in [0,\infty)\}.$$

In what follows, the superscript $n$ is often suppressed for notational simplicity. A (nonlinear) operator $\Delta : \mathbf{L}_{1e} \to \mathbf{L}_{1e}$ is said to be *causal* if $P_T \Delta P_T = P_T \Delta$ for all $T \geq 0$. A

This work was supported by the Swedish Research Council through the LCCC Linnaeus centre.
S. Z. Khong and A. Rantzer are with Department of Automatic Control, Lund University, SE-221 00 Lund, Sweden. seizhen@control.lth.se; rantzer@control.lth.se
C. Briat is with Department of Biosystems Science and Engineering, Swiss Federal Institute of Technology Zürich (ETH Zürich), CH-4058, Basel, Switzerland. corentin.briat@bsse.ethz.ch

causal $\Delta$ is called *bounded* if the Lipschitz bound

$$\|\Delta\| := \sup_{T>0; \|P_T x\|_1 \neq 0} \frac{\|P_T \Delta x\|_1}{\|P_T x\|_1} < \infty. \quad (1)$$

$\Delta$ is said to be *monotone* [1] if for every $u_1, u_2 \in \mathbf{L}_{1e}$,

$$u_1 \geq u_2 \implies \Delta u_1 \geq \Delta u_2,$$

where the inequality is taken coordinate-wise, that is $u_1 \geq u_2$ means each coordinate of $u_1$ is larger than or equal to the corresponding coordinate of $u_2$. In the case where $\Delta$ is linear, monotonicity is equivalent to $\Delta \mathbf{L}_{1e+} \subset \mathbf{L}_{1e+}$. Note that a monotone nonlinear operator $\Delta$ which satisfies $\Delta 0 = 0$ also satisfies $\Delta \mathbf{L}_{1e+} \subset \mathbf{L}_{1e+}$. Such a system is referred to as a *positive* system in this paper. A positive $\Delta$ is said to be *bounded on* $\mathbf{L}_{1e+}$ if

$$\|\Delta\|_+ := \sup_{x \in \mathbf{L}_{1e+}; T>0; \|P_T x\|_1 \neq 0} \frac{\|P_T \Delta x\|_1}{\|P_T x\|_1} < \infty \quad (2)$$

The space of bounded linear causal operators mapping from $\mathbf{L}_1$ to itself is denoted by $\mathscr{L}(\mathbf{L}_1, \mathbf{L}_1)$. Each operator in $\mathscr{L}(\mathbf{L}_1, \mathbf{L}_1)$ has a natural causal extension to an operator mapping from $\mathbf{L}_{1e}$ to $\mathbf{L}_{1e}$, where its $\mathbf{L}_1$-to-$\mathbf{L}_1$ induced norm is equal to the Lipschitz bound of the extension [9, Section 2.4]. Note that $G \in \mathscr{L}(\mathbf{L}_1, \mathbf{L}_1)$ is positive if, and only if, its causal extension satisfies $G\mathbf{L}_{1e+} \subset \mathbf{L}_{1e+}$.

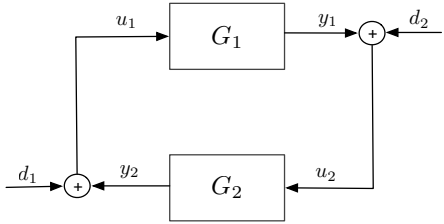

Fig. 1. Positive feedback interconnection of positive systems.

Consider the positive feedback interconnection of nonlinear causal $G_1 : \mathbf{L}_{1e} \to \mathbf{L}_{1e}$ and $G_2 : \mathbf{L}_{1e} \to \mathbf{L}_{1e}$ illustrated in Figure 1. Mathematically,

$$\begin{cases} u_2 = G_1 u_1 + d_2 \\ u_1 = G_2 u_2 + d_1. \end{cases} \quad (3)$$

Denote the feedback interconnection of $G_1$ and $G_2$ by $[G_1, G_2]$,

*Definition 2.1:* $[G_1, G_2]$ is said to be *well-posed* if the map $(u_1, u_2) \mapsto (d_1, d_2)$ defined by (3) has a causal inverse $H$ on $\mathbf{L}_{1e}$. It is *positive* it is well-posed and the inverse $H$ is positive. It is *stable on* $\mathbf{L}_{1e+}$ if it is well-posed, positive, and the inverse $H$ is bounded on $\mathbf{L}_{1e+}$ (cf. (2)). It is *stable* if it is well-posed and the inverse $H$ is bounded (cf. (1)).

Note that if $G_1$ and $G_2$ are bounded and $G_2$ is linear as is the case for standard IQC analysis [8], stability of $[G_1, G_2]$ is equivalent to $(I - G_2 G_1)^{-1}$ being bounded and causal.

*Remark 2.2:* The notion of *stability on* $\mathbf{L}_{1e+}$ is introduced here to characterise the type of well-posed feedback interconnections of positive systems (cf. Figure 1) that remain stable under disturbances $d_1, d_2 \in \mathbf{L}_{1e+}$. In other words, feedback stability is achieved as far as positive signals are concerned. The next lemma shows that for linear systems, this is equivalent to stability on $\mathbf{L}_{1e}$.

*Lemma 2.3:* A linear causal positive $\Delta : \mathbf{L}_{1e} \to \mathbf{L}_{1e}$ is bounded on $\mathbf{L}_{1e+}$ if, and only if, it is bounded on $\mathbf{L}_{1e}$.

*Proof:* Sufficiency is trivial. For necessity, note that any $x \in \mathbf{L}_{1e}$ can be written as $x = x_1 - x_2$, where $x_1, x_2 \in \mathbf{L}_{1e+}$ satisfy $\|P_T x_1\|_1 \leq \|P_T x\|_1$ and $\|P_T x_2\|_1 \leq \|P_T x\|_1$ for all $T > 0$. Moreover, since $\Delta$ is linear, $\Delta x = \Delta x_1 - \Delta x_2$, whereby

$$\begin{aligned} \|P_T \Delta x\|_1 &= \|P_T \Delta P_T x\|_1 \\ &\leq \|P_T \Delta P_T x_1\|_1 + \|P_T \Delta P_T x_2\|_1 \\ &\leq \|\Delta\|_+ \|P_T x_1\|_1 + \|\Delta\|_+ \|P_T x_2\|_1 \\ &\leq 2\|\Delta\|_+ \|P_T x\|_1. \end{aligned}$$

This establishes the boundedness of $\Delta$ on $\mathbf{L}_{1e}$. ∎

*B. Feedback interconnections of positive LTI systems*

Given a matrix $A \in \mathbb{R}^{p \times m}$, $A^T \in \mathbb{R}^{m \times p}$ denotes its transpose. The inequality $A > 0$ ($A \geq 0$) means that all elements of the matrix $X$ are positive (nonnegative). A matrix $A \in \mathbb{R}^{n \times n}$ is called *Hurwitz* if all its eigenvalues have negative real parts. It is said to be *Metzler* if its off-diagonal entries are all nonnegative, i.e. $A_{ij} \geq 0$ for $i \neq j$. The spectrum of a square matrix $M \in \mathbb{C}^{n \times n}$ and spectral radius are defined, respectively, by

$$\mathrm{spec}(M) := \{\lambda \in \mathbb{C} : \lambda I - M \text{ is a singular matrix}\};$$
$$\rho(M) := \max\{|\lambda| : \lambda \in \mathrm{spec}(M)\}.$$

*Proposition 2.4 ([5, Prop. 1 & 2]):* Given a Metzler $A \in \mathbb{R}^{n \times n}$, the following statements are equivalent:

(i) $A$ is Hurwitz;
(ii) There exists a $z \in \mathbb{R}^n$ such that $z > 0$ and $z^T A < 0$;
(iii) $-A^{-1}$ exists and has nonnegative entries.

Given a $B \in \mathbb{R}_+^{n \times n}$, the following statements are equivalent:

(iv) $\rho(B) < 1$;
(v) There exists a $z \in \mathbb{R}^n$ such that $z > 0$ and $z^T B < z^T$;
(vi) $(I - B)^{-1}$ exists and has nonnegative entries.

Some preliminary results on the feedback interconnection of two LTI systems are provided below. They demonstrate the fact that static-gain conditions play a crucial role in the robust feedback stability of positive LTI systems. Suppose that $G_1, G_2$ have real-rational proper realisations:

$$\begin{aligned} \hat{G}_1(s) &= C_1(sI - A_1)^{-1}B_1 + D_1 \quad \text{and} \\ \hat{G}_2(s) &= C_2(sI - A_2)^{-1}B_2 + D_2, \end{aligned} \quad (4)$$

where $A_1$ and $A_2$ are Metzler and $B_1 \geq 0$, $B_2 \geq 0$, $C_1 \geq 0$, $C_2 \geq 0$, $D_1 \geq 0$, and $D_2 \geq 0$. These imply that $G_1$ and $G_2$ are positive [5]. In the following the $\hat{\cdot}$ notation will be used to denote the equivalent frequency-domain representation of an LTI system or the Laplace transform of a signal in $\mathbf{L}_1$ as above.

*Proposition 2.5 ([10, Thm. 4]):* Given $G_1$ and $G_2$ with realisations as in (4), if $D_1 D_2 - I$ (or equivalently, $D_2 D_1 - I$) is Hurwtiz, then $[G_1, G_2]$ is positive. Under this condition,

$[G_1, G_2]$ is stable if, and only if, $A_1$, $A_2$, and $\hat{G}_1(0)\hat{G}_2(0)-I$ (or equivalently, $\hat{G}_2(0)\hat{G}_1(0) - I$) are all Hurwitz.

*Corollary 2.6:* Given $G_1$ and $G_2$ with realisations as in (4), if $\rho(D_1 D_2) < 1$, then $[G_1, G_2]$ is positive. Under this condition, $[G_1, G_2]$ is stable if, and only if, $A_1$ and $A_2$ are Hurwitz and $\rho(\hat{G}_1(0)\hat{G}_2(0)) < 1$.

*Proof:* First note that repeated applications of Proposition 2.4 yields that the Hurwitzness of $D_1 D_2 - I$ is equivalent to $(I - D_1 D_2)^{-1} \geq 0$, which in turn is equivalent to $\rho(D_1 D_2) < 1$. Also by Proposition 2.4(iii), $\hat{G}_1(0)\hat{G}_2(0) \geq 0$ whenever $A_1$ and $A_2$ are Hurwitz. By the same token as before, the Hurwitzness of $\hat{G}_1(0)\hat{G}_2(0) - I$ is equivalent to $\rho(\hat{G}_1(0)\hat{G}_2(0)) < 1$. The claim thus follows from Proposition 2.5. ∎

Denote by $\mathbf{H}_\infty$ the Hardy space of holomorphic transfer functions in the right-half complex plane, which admit equivalent LTI operators mapping from the space of square-integrable functions

$$\mathbf{L}_2 := \left\{ v : [0, \infty) \to \infty : \int_0^\infty |v(t)|_2 \, dt < \infty \right\}$$

to itself via the Laplace transform isomorphism [11, Thm. 3.30], where $|\cdot|_2$ denotes the Euclidean norm.

*Proposition 2.7 ([12, Thm. 3]):* The feedback interconnection of positive $G_1, G_2 \in \mathbf{H}_\infty$ is such that $(I - \hat{G}_2\hat{G}_1)^{-1} \in \mathbf{H}_\infty$ and $(I - G_2 G_1)^{-1}$ is positive on $\mathbf{L}_2$, if and only if $\rho(\hat{G}_1(0)\hat{G}_2(0)) < 1$.

*Remark 2.8:* While Proposition 2.7 is different from the setting of this paper in that the former considers the $\mathbf{L}_2$ space as opposed to $\mathbf{L}_1$, it serves to illustrate the important fact that the robust stability of a positive closed-loop is determined by the static gains of the open-loop systems. Proposition 2.7 can also be seen as a generalisation of aspects of Corollary 2.6 to distributed-parameter transfer functions in $\mathbf{H}_\infty$, which do not admit realisations of the form $C(sI - A)^{-1}B + D$.

In Section IV, it will be demonstrated that the condition $\rho(\hat{G}_1(0)\hat{G}_2(0)) < 1$ can be generalised via integral linear constraints for guaranteeing robust closed-loop stability of positive systems that are time-varying or nonlinear.

## III. POSITIVITY OF NONLINEAR FEEDBACK SYSTEMS

This section demonstrates that a positive feedback interconnection of two positive causal nonlinear (not necessarily bounded) systems can be guaranteed to be positive by the same condition by Willems [9] that ensures the feedback is well-posed. Consider the feedback interconnection of $G_1 : \mathbf{L}_{1e} \to \mathbf{L}_{1e}$ and $G_2 : \mathbf{L}_{1e} \to \mathbf{L}_{1e}$ in Figure 1. It is assumed that for $i = 1, 2$,

(i) $G_i$ is causal and positive;
(ii) $G_i$ is locally Lipschitz continuous on $\mathbf{L}_{1e}$, i.e.

$$\sup_{x, y \in \mathbf{L}_{1e}; P_T x \neq P_T y} \frac{\|P_T(G_i x - G_i y)\|_1}{\|P_T(x - y)\|_1} < \infty, \quad \forall T > 0.$$

The instantaneous gain of $G_i$ is defined by

$$\alpha(G_i) := \sup_{T > 0} \inf_{\Delta T > 0} \sup_{\substack{x, y \in \mathbf{L}_{1e}; P_T x = P_T y \\ P_{T + \Delta T}(x - y) \neq 0}} \frac{\|P_{T + \Delta T}(G_i x - G_i y)\|_1}{\|P_{T + \Delta T}(x - y)\|_1}.$$

Note that $\alpha(G) = 0$ if $G$ is strictly causal, i.e. the input has no instantaneous effect on the output. It is shown in [9, Thm. 4.1] that if $\alpha(G_1)\alpha(G_2) < 1$, then the $[G_1, G_2]$ is well-posed as per Definition 2.1. The following demonstrates that $\alpha(G_1)\alpha(G_2) < 1$ also guarantees positivity of $[G_1, G_2]$.

*Theorem 3.1:* If $\alpha(G_1)\alpha(G_2) < 1$, then $[G_1, G_2]$ is well-posed and positive.

*Proof:* Well-posedness is established in [9, Thm. 4.1]. As for positivity, let $u := (u_1, u_2)$, $d := (d_1, d_2)$, and $G(u_1, u_2) := (G_2 u_2, G_1 u_1)$. The feedback equation (3) can be written as

$$u = Gu + d.$$

Well-posedness and positivity of $[G_1, G_2]$ is thus equivalent to the causal positive invertibility of $(I - G)$. [9, Thm. 4.1] establishes that when $\alpha(G_1)\alpha(G_2) < 1$, $P_{T_0 + \Delta T} G P_{T_0 + \Delta T}$ is contractive for $T_0 \geq 0$ and sufficiently small $\Delta T > 0$, whereby $P_{T_0 + \Delta T}(I - G)P_{T_0 + \Delta T}$ is invertible. Invertibility of $I - G$ is then shown on consecutive intervals that cover the whole half line $[0, \infty)$. Moreover, causality of $(I - G)^{-1}$ is established from the fact that $G$ is causal and the invertibility of $I - G$ is based on the convergence of successive approximations

$$u_{n+1} = Gu_n + d \quad n = 0, 1, \ldots$$

for arbitrary $u_0$. Therefore, when $G$ is positive, $d \geq 0$, and $u_0 \geq 0$, it follows that $\lim_{n \to \infty} u_n \geq 0$. This establishes positivity of $(I - G)^{-1}$ along the same lines of causality. ∎

Consider again the setting of Proposition 2.5, where $\hat{G}_1(s) = C_1(sI - A_1)^{-1} B_1 + D_1$ and $\hat{G}_2(s) = C_2(sI - A_2)^{-1} B_2 + D_2$ with $A_1$ and $A_2$ Metzler and $B_1 \geq 0$, $B_2 \geq 0$, $C_1 \geq 0$, $C_2 \geq 0$, $D_1 \geq 0$, and $D_2 \geq 0$. Observe that

$$\alpha(G_i) = \sup_{|u|=1} |D_i u| =: \|D_i\|_{1 \to 1}.$$

Thus, Theorem 3.1 ensures that $[G_1, G_2]$ is positive whenever $\|D_1\|_{1 \to 1} \|D_2\|_{1 \to 1} < 1$. Alternatively, it may be established by the first part of Corollary 2.6. To see this, note that

$$\rho(D_1 D_2) \leq \|D_1 D_2\|_{1 \to 1} \leq \|D_1\|_{1 \to 1} \|D_2\|_{1 \to 1} < 1.$$

Theorem 3.1 is useful in that it can be applied to nonlinear time-varying systems.

## IV. INTEGRAL LINEAR CONSTRAINTS

Given a causal bounded $\Delta : \mathbf{L}_{1e} \to \mathbf{L}_{1e}$, define the graph of $\Delta$ with respect to positive signals as

$$\mathcal{G}(\Delta) := \left\{ \begin{bmatrix} x \\ y \end{bmatrix} \in \mathbf{L}_{1+} : y = \Delta x \right\}.$$

Similarly, define the inverse graph as

$$\mathcal{G}'(\Delta) := \left\{ \begin{bmatrix} y \\ x \end{bmatrix} \in \mathbf{L}_{1+} : x = \Delta y \right\}.$$

The following definition of integral linear constraint follows that in [3]. As will be seen later, this allows the establishment of a link from the main robust stability result to the static gain condition in Corollary 2.6.

*Definition 4.1:* A causal bounded system $\Delta : \mathbf{L}_{1e}^m \to \mathbf{L}_{1e}^p$ is said to satisfy the integral linear constraint (ILC) defined by the multiplier $\Pi \in \mathbb{R}^{1 \times (m+p)}$ if

$$\int_0^\infty \Pi v(t)\, dt \geq 0 \quad \forall v \in \mathcal{G}(\Delta).$$

This is denoted $\Delta \in \text{ILC}(\Pi)$. On the contrary, $\Delta$ is said to satisfy the strict complementary ILC if

$$\int_0^\infty \Pi w(t)\, dt \leq -\epsilon \int_0^\infty 1_{p+m}^T w(t)\, dt \quad \forall w \in \mathcal{G}'(\Delta)$$

for some $\epsilon > 0$, where $1_n \in \mathbb{R}^n$ denotes the vector whose entries are all ones. This is denoted $\Delta \in \text{ILC}^c(\Pi)$.

*Remark 4.2:* Having $\Pi$ as a bounded LTI operator from $\mathbf{L}_1$ to $\mathbf{L}_1$ does not generalise the definition of an ILC for the following reason. Denote by $\pi$ the convolution kernel corresponding to $\Pi$ and $*$ the convolution operation, then

$$\int_0^\infty \left( \pi * \begin{bmatrix} u \\ \Delta u \end{bmatrix} \right)(t)\, dt$$
$$= \int_0^\infty \left( \pi * \begin{bmatrix} u \\ \Delta u \end{bmatrix} \right) e^{-st}\, dt \bigg|_{s=0}$$
$$= \hat{\pi}(0) \begin{bmatrix} \hat{u}(0) \\ \widehat{\Delta u}(0) \end{bmatrix}.$$

That is, only the static gain of $\hat{\pi}$, which is a constant matrix, matters in the value of the integral.

*Lemma 4.3:* Suppose $\Delta : \mathbf{L}_{1e}^m \to \mathbf{L}_{1e}^p$ is a causal bounded positive LTI system, then $\Delta \in \text{ILC}(\Pi)$ is equivalent to

$$\Pi \begin{bmatrix} I_m \\ \hat{\Delta}(0) \end{bmatrix} \geq 0, \tag{5}$$

where $\hat{\Delta}$ denotes the transfer function representation of $\Delta$. Similarly, $\Delta \in \text{ILC}^c(\Pi)$ is equivalent to

$$\Pi \begin{bmatrix} \hat{\Delta}(0) \\ I_m \end{bmatrix} < 0.$$

*Proof:* For any $u \in \mathbf{L}_{1+}$,

$$\int_0^\infty \Pi \begin{bmatrix} u(t) \\ (\Delta u)(t) \end{bmatrix} dt$$
$$= \int_0^\infty \Pi \begin{bmatrix} u(t) \\ (\Delta u)(t) \end{bmatrix} e^{-st}\, dt \bigg|_{s=0}$$
$$= \Pi \begin{bmatrix} I_m \\ \hat{\Delta}(0) \end{bmatrix} \hat{u}(0).$$

Now note that $\Pi \begin{bmatrix} I_m \\ \hat{\Delta}(0) \end{bmatrix} \hat{u}(0) \geq 0$ for all $u \in \mathbf{L}_{1+}$ is equivalent to $\Pi \begin{bmatrix} I_m \\ \hat{\Delta}(0) \end{bmatrix} \geq 0$, since $\hat{u}(0) = \int_0^\infty u(t)\, dt \geq 0$. The second part of the claim can be shown using the same lines of arguments. ∎

Note that the static-gain condition (5) in Lemma 4.3 is amenable to distributed verifications, since it only involves multiple vector multiplications. The computation also scales linearly with the dimensions of the system concerned, much in the spirit of [5].

## V. ROBUSTNESS ANALYSIS

The main robust stability result for positive feedback systems is stated. Note that the open-loop systems do *not* need to be positive for the following to hold.

*Theorem 5.1:* Given a (nonlinear) bounded causal $G_1 : \mathbf{L}_{1e}^m \to \mathbf{L}_{1e}^p$ and a linear bounded causal $G_2 : \mathbf{L}_{1e}^p \to \mathbf{L}_{1e}^m$, suppose there exists a $\Pi \in \mathbb{R}^{1 \times m+p}$ such that
  (i) $[\tau G_1, G_2]$ is well-posed and positive for all $\tau \in [0,1]$;
  (ii) $\tau G_1 \in \text{ILC}(\Pi)$ for all $\tau \in [0,1]$; and
  (iii) $G_2 \in \text{ILC}^c(\Pi)$.
Then $[G_1, G_2]$ is stable on $\mathbf{L}_{1e+}$. In addition, if $G_1$ is also linear, then $[G_1, G_2]$ is stable.

*Proof:* By hypothesis, for all $\tau \in [0,1]$, $v \in \mathcal{G}(\tau G_1)$ and $w \in \mathcal{G}'(G_2)$, there exists $\epsilon > 0$ such that

$$\int_0^\infty \Pi v(t)\, dt \geq 0$$

and

$$\int_0^\infty \Pi w(t)\, dt \leq -\epsilon \int_0^\infty 1_{p+m}^T w(t)\, dt.$$

Define $\Psi := 2\Pi + \epsilon 1_{p+m}^T$. The inequalities can be written as

$$\int_0^\infty \Psi v(t)\, dt \geq \epsilon \int_0^\infty 1_{p+m}^T v(t)\, dt = \epsilon \|v\|_1$$

and

$$\int_0^\infty \Psi w(t)\, dt \leq -\epsilon \int_0^\infty 1_{p+m}^T w(t)\, dt = -\epsilon \|w\|_1,$$

where the fact that $v, w \in \mathbf{L}_{1+}$ has been exploited. Therefore, by the feedback configuration in (3),

$$\epsilon(\|v\|_1 + \|w\|_1) \leq \int_0^\infty \Psi(v(t) - w(t))\, dt$$
$$= \int_0^\infty \Psi \begin{bmatrix} d_1(t) \\ -d_2(t) \end{bmatrix} dt$$
$$\leq \int_0^\infty \bar{\psi} 1_{p+m}^T \begin{bmatrix} d_1(t) \\ d_2(t) \end{bmatrix} dt$$
$$= \bar{\psi}(\|d_1\|_1 + \|d_2\|_1),$$

for all $v \in \mathcal{G}(\tau G_1)$, $w \in \mathcal{G}'(G_2)$, and $d_1, d_2 \in \mathbf{L}_{1+}$, where $\bar{\psi} := \max_{i=1,\ldots,m+p} |\psi_i|$ and $\Psi = [\psi_1, \ldots \psi_{m+p}]$. In other words,

$$\|v\|_1 + \|w\|_1 \leq \frac{\bar{\psi}}{\epsilon}(\|d_1\|_1 + \|d_2\|_1). \tag{6}$$

Note that when $d_2 = 0$, $d_1 = (I - \tau G_2 G_1)u_1$ by linearity of $G_2$. As such, it follows from (6) and $v = \begin{bmatrix} u_1 \\ \tau G_1 u_1 \end{bmatrix}$ that for all $u_1 \in \mathbf{L}_{1+}$,

$$\|u_1\| \leq \frac{\bar{\psi}}{\epsilon} \|(I - \tau G_2 G_1)u_1\|_1. \tag{7}$$

By the well-posedness and positivity assumption, the inverse $(I - \tau G_2 G_1)^{-1}$ is well-defined on $\mathbf{L}_{1e}$ and satisfies $(I - \tau G_2 G_1)^{-1} \mathbf{L}_{1e+} \subset \mathbf{L}_{1e+}$. Given any $u \in \mathbf{L}_{1e+}$ and $\nu \in [0,1]$, define

$$u_T := (I - \nu G_2 G_1)^{-1} P_T (I - \nu G_2 G_1) u \in \mathbf{L}_{1e+}.$$

Then
$$\begin{aligned}\|P_T u\|_1 &= \|P_T u_T\|_1 \\ &\leq \|u_T\|_1 \\ &\leq \frac{\bar\psi}{\epsilon}\|(I-\nu G_2 G_1)u_T\|_1 \\ &= \frac{\bar\psi}{\epsilon}\|P_T(I-\nu G_2 G_1)u\|_1 \\ &= \frac{\bar\psi}{\epsilon}\|P_T(I-\tau G_2 G_1)u + (\tau-\nu)P_T G_2 G_1 u\|_1 \\ &\leq \frac{\bar\psi}{\epsilon}\|P_T(I-\tau G_2 G_1)u\|_1 \\ &\quad + \frac{\bar\psi}{\epsilon}|\tau-\nu|\|G_2\|\|G_1\|\|P_T u\|_1,\end{aligned} \quad (8)$$

where (7) has been used to arrive at the second inequality. Since $(I-\nu G_2 G_1)^{-1}$ is bounded on $\mathbf{L}_{1e+}$ for $\nu=0$, it follows from (8) that $(I-\tau G_2 G_1)^{-1}$ is bounded on $\mathbf{L}_{1e+}$ for all $\tau < (\frac{\bar\psi}{\epsilon}\|G_2\|\|G_1\|)^{-1}$. By an inductive argument and repeated applications of (8), it follows that $(I-\tau G_2 G_1)^{-1}$ is bounded on $\mathbf{L}_{1e+}$ for all $\tau \in [0,1]$. In other words, $[\tau G_1, G_2]$ is stable on $\mathbf{L}_{1e+}$ for all $\tau \in [0,1]$.

For the final claim, note that linearity of $G_1$ and $G_2$ implies that of $(I - G_2 G_1)^{-1}$ and apply Lemma 2.3. ∎

## VI. EXAMPLES

### A. LTI systems

Consider the feedback interconnection of two positive LTI systems $G_1$ and $G_2$, where $\hat{G}_1(s) = C_1(sI-A_1)^{-1}B_1 + D_1$ and $\hat{G}_2(s) = C_2(sI-A_2)^{-1}B_2 + D_2$ with $A_1$ and $A_2$ Metzler and $B_1 \geq 0$, $B_2 \geq 0$, $C_1 \geq 0$, $C_2 \geq 0$, $D_1 \geq 0$, and $D_2 \geq 0$. Suppose that $\|D_1\|_{1\to 1}\|D_2\|_{1\to 1} < 1$, whereby

$$\|\tau D_1\|_{1\to 1}\|D_2\|_{1\to 1} < 1 \;\forall \tau \in [0,1],$$

so that $[\tau G_1, G_2]$ is well-posed and positive by Theorem 3.1. It is shown below that the sufficiency part of Corollary 2.6, i.e. $\rho(\hat{G}_1(0)\hat{G}_2(0)) < 1$ implies that $[G_1, G_2]$ is stable, can be recovered from the main result Theorem 5.1 using a particular multiplier $\Pi$.

To begin with, note that $\rho(\hat{G}_1(0)\hat{G}_2(0)) < 1$ is equivalent to the existence of a $z > 0$ such that $z^T(\hat{G}_1(0)\hat{G}_2(0) - I) < 0$ by Proposition 2.4. Define

$$\Pi := z^T \begin{bmatrix} \hat{G}_1(0) & -I \end{bmatrix}.$$

It is straightforward to see that

$$\Pi \begin{bmatrix} \hat{G}_2(0) \\ I \end{bmatrix} < 0 \quad \text{and} \quad \Pi \begin{bmatrix} I \\ \tau \hat{G}_1(0) \end{bmatrix} \geq 0 \;\forall \tau \in [0,1].$$

By Lemma 4.3, these are equivalent to $\tau G_1 \in \mathrm{ILC}(\Pi)$ for all $\tau \in [0,1]$ and $G_2 \in \mathrm{ILC}^c(\Pi)$. As such, $[G_1, G_2]$ is stable by Theorem 5.1.

### B. Feedback channels with time-varying gain

Suppose two strictly causal bounded positive LTI systems $M_1$ and $M_2$ are connected in a positive feedback loop via channels with time-varying gains $\delta_1 I$ and $\delta_2 I$ satisfying $0 \leq \delta_1(t) \leq \alpha$ and $0 \leq \delta_2(t) \leq \beta$. Stability of this is equivalent to the stability of the structured feedback interconnection $[M,\Delta]$, where

$$M = \begin{bmatrix} 0 & M_1 \\ M_2 & 0 \end{bmatrix} \quad \text{and} \quad \Delta = \begin{bmatrix} \delta_1 I & 0 \\ 0 & \delta_2 I \end{bmatrix}.$$

Note that with $\Pi := (\alpha I, \beta I, -I, -I)$, $\tau\Delta \in \mathrm{ILC}(\Pi) \,\forall \tau \in [0,1]$. By the strict causality assumption on $M$, application of Theorem 3.1 yields that $[M, \tau\Delta]$ is well-posed and positive for all $\tau \in [0,1]$. It thus follows from Lemma 4.3 and Theorem 5.1 that $[M,\Delta]$ is stable if

$$\Pi \begin{bmatrix} \hat{M}(0) \\ I \end{bmatrix} < 0.$$

## VII. CONCLUSIONS

The positivity of positive feedback interconnections of positive nonlinear systems is shown to be guaranteed by an instantaneous gain condition. The static gain conditions on robust feedback stability of positive LTI systems have been generalised to time-varying and nonlinear systems via the notion of integral linear constraints. Future research directions may involve investigating if the use of linear time-varying multipliers reduces conservatism in robustness analysis and the use of $\mathbf{L}_\infty$ based conditions as a replacement of the $\mathbf{L}_1$ based ILC.

## VIII. ACKNOWLEDGEMENT

The first author is grateful to Chung-Yao Kao for technical discussions about the paper.

## REFERENCES

[1] D. Angeli and E. D. Sontag, "Monotone control systems," *IEEE Trans. Autom. Contr.*, vol. 48, no. 10, pp. 1684–1698, 2003.
[2] M. A. Rami, "Solvability of static output-feedback stabilization for LTI positive systems," *Systems and Control Letters*, vol. 60, pp. 704–708, 2011.
[3] C. Briat, "Robust stability and stabilization of uncertain linear positive systems via integral linear constraints: $L_1$-gain and $L_\infty$-gain characterization," *International Journal of Robust and Nonlinear Control*, vol. 23, pp. 1932–1954, 2013.
[4] T. Tanaka and C. Langbort, "The bounded real lemma for internally positive systems and H-infinity structured static state feedback." *IEEE Trans. Autom. Contr.*, vol. 56, no. 9, pp. 2218–2223, 2011.
[5] A. Rantzer, "Distributed control of positive systems," *ArXiv e-prints*, 2014. [Online]. Available: http://arxiv.org/abs/1203.0047
[6] M. A. Rami, *Positive Systems*, ser. Lecture Notes in Control and Information Sciences. Springer, 2009, vol. 389, ch. Stability Analysis and Synthesis for Linear Positive Systems with Time-Varying Delays, pp. 205–215.
[7] F. Forni and R. Sepulchre, "Differentiall positive systems," *IEEE Trans. Autom. Contr.*, 2015, in press.
[8] A. Megretski and A. Rantzer, "System analysis via integral quadratic constraints," *IEEE Trans. Autom. Contr.*, vol. 42, no. 6, pp. 819–830, 1997.
[9] J. C. Willems, *The Analysis of Feedback Systems*. Cambridge, Massachusetts: MIT Press, 1971.
[10] Y. Ebihara, D. Peaucelle, and D. Arzelier, "$L_1$ gain analysis of linear positive systems and its application," in *IEEE Conference on Decision and Control*, 2011, pp. 4029–4034.
[11] G. E. Dullerud and F. Paganini, *A Course in Robust Control Theory: A Convex Approach*, ser. Texts in Applied Mathematics 36. Springer, 2000.
[12] T. Tanaka, C. Langbort, and V. Ugrinovskii, "DC-dominant property of cone-preserving transfer functions," *Systems and Control Letters*, vol. 62, pp. 699–707, 2013.